\documentstyle[12pt]{article}

\begin{document}

\title{ \Large\bf A Remark on Generalized Covering Groups}
\author{ Behrooz Mashayekhy\\ Department of Mathematics\\ Ferdowsi University of
Mashhad \\ P.O.Box 1159-91775 \\ Mashhad, Iran\\ E-mail:
mashaf@science2.um.ac.ir }
\date{}
\maketitle
\begin{abstract}
 Let ${\cal N}_c$ be the variety of nilpotent groups of class at most $c\ \
(c\geq 2)$ and $G=Z_r\oplus Z_s $ be the direct sum of two finite
cyclic groups. It is shown that if the greatest common divisor of
$r$ and $s$ is not one, then $G$ does not have any ${\cal
N}_c$-covering group for every $c\geq 2$.  This result gives an
idea that Lemma 2 of J.Wiegold [6] and Haebich's Theorem [1], a
vast generalization of the Wiegold's Theorem, can {\it not} be
generalized to the variety of nilpotent groups of class at most
$c\geq 2$.\\
A.M.S.Classification (1990): 20F12,20F18,20K25 
\end{abstract}
\begin{center} {\bf 1. Notation and Preliminaries }\end{center}

 We assume that the reader is familiar with the notiona of the verbal subgroup,
 $V(G)$ , and the marginal subgroup, $V^*(G)$ , associated with a variety of
groups, ${\cal V}$, and a group $G$ , and the basic commutators (see [2]).

 Let ${\cal V}$ be a variety of groups, and $G$ be a group with the following
free presentation , $ 1\rightarrow R\rightarrow F\rightarrow G\rightarrow 1 $\
. Then the {\it Baer-invariant} of $G$ , with respect to the variety ${\cal
V}$ , denoted by ${\cal V}M(G)$ , is defined to be  $ R\cap V(F)/[RV^*F] $
where $V(F)$ is the verbal subgroup of $F$ and
\[ [RV^*F]=<v(f_1,\ldots ,f_{i-1},f_ir,f_{i+1},\ldots ,f_n)v(f_1,\ldots
,f_i,\ldots ,f_n)^{-1}\mid \] \[ r\in R\ ,\ f_i\in F,\ v\in V\ ,\ 1\leq i\leq
n\ ,\ n\in {\bf N}>\ .\]
 In particular, if ${\cal V}$ is the variety of nilpotent groups of class at
most $c\ \ (c\geq 1)$ , then
$$ {\cal V}M(G)=\frac {R\cap \gamma_{c+1}(F)} {[R,\ _cF]}\ .$$

 A ${\cal V}$\_stem cover of $G$ is an exact sequence $ 1\rightarrow
A\rightarrow G^*\rightarrow G\rightarrow 1 $ such that $A\subseteq V(G^*)\cap
V^*(G^*)$ ,where $G^*$ is said to be a ${\cal V}$-covering group of $G$ . It
is of interest to know which class of groups does {\it not} have a ${\cal
V}$-covering group. For further details see C.R.Leedham-Green and S.McKay
[4].
\newpage
\begin{center} {\bf 2. The Main Result}  \end{center}

 The following theorem is important in our study.\\
{\bf Theorem 1 }

 Let ${\cal N}_c$ be the variety of nilpotent groups of class at most $c\ \
(c\geq 1)$ and $r,s$ be two positive integers with $(r,s)=d$ their
greatest common divisor. If $G=Z_r\oplus Z_s$, then ${\cal N}_cM(G)\cong
Z_d\oplus Z_d\oplus \ldots \oplus Z_d $ (n-copies) , where $n$ is the number
of basic commutators of weight $c+1$ on two letters.\\
{\bf Proof.}

 Take the following free presentation for $G$ :
$$ 1\longrightarrow R\longrightarrow F\longrightarrow G\longrightarrow 1\ , $$
where $F$ is the free group on $\{x_1,x_2\}$ and
$R=<{x_1}^{r},{x_2}^{s},\gamma_2(F)>\ .$ Clearly $R=S\gamma_2(F)$ , where
$S=<{x_1}^r,{x_2}^s>^F.$ So the Baer-invariant of $G$ with respect to the
variety ${\cal N}_c$ is
$$ {\cal N}_cM(G)=\frac {R\cap \gamma_{c+1}(F)} {[R,\ _cF]} =\frac
{S\gamma_2(F)\cap \gamma_{c+1}(F)} {[s\gamma_2(F),\ _cF]}$$  $$ =\frac
{\gamma_{c+1}(F)}
{[S,\ _cF]\gamma_{c+2}(F)}\cong \frac {\gamma_{c+1}(F)/\gamma_{c+2}(F)}
{[S,\ _cF]\gamma_{c+2}(F)/\gamma_{c+2}(F)}$$

 By P.Hall's Theorem (see[2]) $\gamma_{c+1}(F)/\gamma_{c+2}(F) $ is a free
abelian group freely generated by all the basic commutators of weight $c+1$ on
two letters. For all $a_i\in F$ and any $k\in {\bf Z}$ we have
$$ [a_1,\ldots ,{a_i}^k,\ldots ,a_{c+1}]=[a_1,\ldots ,a_i,\ldots ,a_{c+1}]^k
\ \ \ \ \ ( mod\ \ \gamma_{c+2}(F))\ .$$
 Hence $ [S,_cF]\gamma_{c+2}(F)/\gamma_{c+2}(F) $ is a free abelian group
freely generated by the set of all $d$-th powers of all basic commutators of
weight $c+1$ on two letters. Hence ${\cal N}_cM(G)\cong Z_d\oplus \ldots \oplus
Z_d $ ($n$-copies) where n is the number of basic commutators of weight $c+1$
on two letters.\ $\Box $

 The above theorem is somehow a generalization of M.R.R.Moghaddam's Theorem
(see [5]).

  Now we are in a position to state the main theorem.\\
 {\bf Theorem 2}

 Let ${\cal N}_c$ be the variety of nilpotent groups of class at most $c\ \
(c\geq 2)$ , and let $(r,s)=d\not =1$ be the greatest common divisor of the
positive integers $r,s$. If $G=Z_r\oplus Z_s $ , then $G$ does {\it not} admit
any ${\cal N}_c$-stem cover for every $c\geq 2$ and so $G$ has no ${\cal
N}_c$-covering group.\\
{\bf Proof.}

 Assume by way of cotradiction that $ 1\rightarrow A\rightarrow
G^*\rightarrow G\rightarrow 1 $ is an ${\cal N}_c$-stem cover for $G$ . Thus
$ A\subseteq \gamma_{c+1}(G^*)\cap Z_c(G^*)$ (where $Z_c(G^*)$ is the $c$-th
center of $G^*$) , $G^*/A\cong G$ and $ A\cong {\cal N}_cM(G)\ .$ Since $G$ is
abelian, $\gamma_2(G^*)\subseteq A$. Therefore $\gamma_{c+2}(G^*)=1 $ , using
the fact $\gamma_2(G^*)\subseteq Z_c(G^*) $ . It implies that
$$ 1=\gamma_{c+2}(G^*)=[\gamma_{c+1}(G^*),G^*]\supseteq
[\gamma_2(G^*),G^*]=\gamma_3(G^*)\ .$$
As $c\geq 2$ , we obtain $\gamma_{c+1}(G^*)=1$ and hence ${\cal N}_cM(G)\cong
A=1 $ which is a contradiction by Theorem 1 .\ $\Box $
\newpage
{\bf Remark}

$(i)$ If $c=1$ then one may construct a covering group for $G$ . (See
G.Karpilovsky [3])

$(ii)$ J.Wiegold in [6] showed that if $G_1,G_2$ are two finite
groups and $G^*_1,G^*_2$ any covering groups of $G_1,G_2$ , respectively, then
the second nilpotent product $G^*_1\stackrel{2}{*} G^*_2$ is a covering group
of $G_1\times G_2$ . Also W.Haebich constructed a covering group for a regular
product of a class of groups (see [1]). Now our Theorem 2 shows that these
notions can not be generalized to the variety of nilpotent groups of class
$c\geq 2$. Because, in general, the direct product $G_1\times G_2$ may not have
an ${\cal N}_c$-covering group ($c\geq 2$).

\end{document}